\documentclass[12pt]{article}
\author{Edward Hanson}
\title{\bf A characterization of Leonard pairs\\ using the notion of a tail}
\date{}
\usepackage{amssymb}
\usepackage{amsmath}
\usepackage{cite}

\newtheorem{definition}{Definition}[section]
\newtheorem{theorem}[definition]{Theorem}
\newtheorem{proposition}[definition]{Proposition}
\newtheorem{lemma}[definition]{Lemma}
\newtheorem{corollary}[definition]{Corollary}

\newtheorem{note}[definition]{Note}
\newtheorem{assumption}[definition]{Assumption}

\usepackage[margin=1.0in]{geometry}

\def\fld{\mathbb K}

\begin{document}
\maketitle

\begin{abstract}
\noindent
Let $V$ denote a vector space with finite positive dimension. We consider an ordered pair of linear transformations 
$A: V\rightarrow V$ and $A^*: V\rightarrow V$ that satisfy (i) and (ii) below:
\begin{enumerate}
\item There exists a basis for $V$ with respect to which the matrix representing $A$ is irreducible tridiagonal and the matrix representing $A^*$ is diagonal.
\item There exists a basis for $V$ with respect to which the matrix representing $A^*$ is irreducible tridiagonal and the matrix representing $A$ is diagonal.
\end{enumerate}
We call such a pair a {\it Leonard pair} on $V$. In this paper, we characterize the Leonard pairs using the notion of a tail. This notion is borrowed from algebraic graph theory.

\end{abstract}

\section{Introduction} \label{sec:intro}

We begin by recalling the notion of a Leonard pair \cite{T:subconst1, T:Leonard, T:qSerre, T:24points, T:intro, T:TD-D}. We will use the following terms. Let $X$ denote a square matrix. Then $X$ is called {\it tridiagonal} whenever each nonzero entry lies on either the diagonal, the subdiagonal, or the superdiagonal. Assume $X$ is tridiagonal. Then $X$ is called {\it irreducible} whenever each entry on
the subdiagonal is nonzero and each entry on the superdiagonal is nonzero.

\medskip
\noindent 
We now define a Leonard pair. For the rest of this paper, $\fld$ will denote a field.

\begin{definition} \label{def:lp} \rm \cite[Definition 1.1]{T:Leonard}
Let $V$ denote a vector space over $\fld$ with finite positive dimension. By a {\it Leonard pair} on $V$, we mean an ordered pair of linear transformations 
$A: V\rightarrow V$ and $A^*: V\rightarrow V$ that satisfy (i) and (ii) below:
\begin{enumerate}
\item There exists a basis for $V$ with respect to which the matrix representing $A$ is irreducible tridiagonal and the matrix representing $A^*$ is diagonal.
\item There exists a basis for $V$ with respect to which the matrix representing $A^*$ is irreducible tridiagonal and the matrix representing $A$ is diagonal.
\end{enumerate}
\end{definition}

\begin{note}
\rm
It is a common notational convention to use $A^{*}$ to represent the conjugate-transpose of $A$. We are not using this convention. In a Leonard pair $A,A^{*}$, the linear transformations $A$ and $A^*$ are arbitrary subject to (i), (ii) above.
\end{note}

\noindent In this paper, we will characterize the Leonard pairs using the notion of a tail. This notion is from algebraic graph theory or, more precisely,
the theory of distance-regular graphs \cite{BIbook, BCN}. The notion was introduced by M.S. Lang \cite{Lang} and developed further in \cite{JTZ}. Our main result, which is Theorem \ref{thm:main} below, can be viewed as an algebraic version of \cite[Theorem 1.1]{JTZ}.

\section{Leonard systems} \label{sec:LS}

\medskip
\noindent When working with a Leonard pair, it is often convenient to consider a closely related object called a Leonard system. To prepare for our definition of a Leonard system, we recall a few concepts from linear algebra. From now on, we fix a nonnegative integer $d$. Let $\mbox{Mat}_{d+1}(\fld)$ denote the $\fld$-algebra consisting of all $d+1$ by $d+1$ matrices with entries in $\fld$. We index the rows and columns by $0,1,\ldots ,d$. We let $\fld^{d+1}$ denote the $\fld$-vector space consisting of all $d+1$ by $1$ matrices with entries in $\fld$. We index the rows by $0,1,\ldots ,d$. Recall that $\mbox{Mat}_{d+1}(\fld)$ acts on $\fld^{d+1}$ by left multiplication. Let $V$ denote a vector space over $\fld$ with dimension $d+1$. Let $\mbox{End}(V)$ denote the $\fld$-algebra consisting of all linear transformations from $V$ to $V$. For convenience, we abbreviate $\mathcal{A}=\mbox{End}(V)$. Observe that $\mathcal{A}$ is $\fld$-algebra isomorphic to $\mbox{Mat}_{d+1}(\fld)$ and that $V$ is irreducible as an $\mathcal{A}$-module. The identity of $\mathcal{A}$ will be denoted by $I$. Let $\{v_{i}\}_{i=0}^{d}$ denote a basis for $V$. For $X\in \mathcal{A}$ and $Y\in \mbox{Mat}_{d+1}(\fld)$, we say that $Y$ {\it represents} $X$ {\it with respect to} $\{v_{i}\}_{i=0}^{d}$ whenever $Xv_{j}=\sum_{i=0}^{d}Y_{ij}v_{i}$ for $0\leq j\leq d$. Let $A$ denote an element of $\mathcal{A}$. A subspace $W\subseteq V$ will be called an {\it eigenspace} of $A$ whenever $W\neq 0$ and there exists $\theta \in \fld$ such that $W=\{v\in V|Av=\theta v\}$; in this case, $\theta$ is the {\it eigenvalue} of $A$ associated with $W$. We say that $A$ is {\it diagonalizable} whenever $V$ is spanned by the eigenspaces of $A$. We say that $A$ is {\it multiplicity-free} whenever it has $d+1$ mutually distinct eigenvalues in $\fld$. Note that if $A$ is multiplicity-free, then $A$ is diagonalizable.

\begin{definition} \label{def:MOidem}
\rm
By a {\it system of mutually orthogonal idempotents} in $\mathcal{A}$, we mean a sequence $\{E_{i}\}_{i=0}^{d}$ of elements in $\mathcal{A}$ such that
\begin{equation}
E_{i}E_{j}=\delta_{i,j}E_{i}\qquad \qquad (0 \leq i,j \leq d), \notag
\end{equation}
\begin{equation}
rank(E_{i})=1\qquad \qquad (0 \leq i \leq d). \notag
\end{equation}
\end{definition}

\begin{definition} \label{def:decomp}
\rm
By a {\it decomposition of $V$}, we mean a sequence $\{U_{i}\}_{i=0}^{d}$ consisting of one-dimensional subspaces of $V$ such that
\begin{equation}
V=\sum_{i=0}^{d}U_{i}\qquad \qquad \text{(direct sum)}. \notag
\end{equation}
\end{definition}

\noindent Definitions \ref{def:MOidem} and \ref{def:decomp} are related in the following lemma, whose proof is left as an exercise.

\begin{lemma} \label{lem:MOidem2}
Let $\{U_{i}\}_{i=0}^{d}$ denote a decomposition of $V$. For $0\leq i\leq d$, define $E_{i}\in \mathcal{A}$ such that $(E_{i}-I)U_{i}=0$ and $E_{i}U_{j}=0$ if $j\ne i$ $(0\leq j\leq d)$. Then $\{E_{i}\}_{i=0}^{d}$ is a system of mutually orthogonal idempotents. Conversely, given a system of mutually orthogonal idempotents $\{E_{i}\}_{i=0}^{d}$ in $\mathcal{A}$, define $U_{i}=E_{i}V$ for $0\leq i\leq d$. Then $\{U_{i}\}_{i=0}^{d}$ is a decomposition of $V$.
\end{lemma}

\begin{lemma} \label{lem:EsumI}
Let $\{E_{i}\}_{i=0}^{d}$ denote a system of mutually orthogonal idempotents in $\mathcal{A}$. Then $I=\sum_{i=0}^{d}E_{i}$.
\end{lemma}
\noindent {\it Proof:} By Lemma \ref{lem:MOidem2}, the sequence $\{E_{j}V\}_{j=0}^{d}$ is a decomposition of $V$. Observe that $\sum_{i=0}^{d}E_{i}$ acts as the identity on $E_{j}V$ for $0\leq j\leq d$. The result follows. \hfill $\Box$ \\

\noindent Let $A$ denote a multiplicity-free element of $\mathcal{A}$ and let $\{\theta_{i}\}^{d}_{i=0}$ denote an ordering of the eigenvalues of $A$. For $0\leq i\leq d$, let $U_{i}$ denote the eigenspace of $A$ for $\theta_{i}$. Then $\{U_{i}\}_{i=0}^{d}$ is a decomposition of $V$; let $\{E_{i}\}_{i=0}^{d}$ denote the corresponding system of idempotents from Lemma \ref{lem:MOidem2}. One checks that $A=\sum_{i=0}^{d}\theta_{i}E_{i}$ and $AE_{i}=E_{i}A=\theta_{i}E_{i}$ for $0\leq i\leq d$. Moreover,
\begin{equation}
E_{i}=\prod_{\genfrac{}{}{0pt}{}{0 \leq  j \leq d}{j\not=i}}\frac{A-\theta_{j}I}{\theta_{i}-\theta_{j}}\qquad \qquad (0\leq i\leq d). \label{eq:EpolyA}
\end{equation}
We refer to $E_{i}$ as the {\it primitive idempotent} of $A$ corresponding to $U_{i}$ (or $\theta_{i}$).

\medskip
\noindent We now define a Leonard system.

\begin{definition} \label{def:ls} \rm \cite[Definition 1.4]{T:Leonard}
By a {\it Leonard system} on $V$, we mean a sequence 
\begin{equation}
(A; \{E_{i}\}_{i=0}^{d}; A^{*}; \{E^{*}_{i}\}_{i=0}^{d}) \notag
\end{equation}
which satisfies (i)--(v) below.
\begin{enumerate}
\item Each of $A,A^{*}$ is a multiplicity-free element of $\mathcal{A}$.
\item $\{E_{i}\}_{i=0}^{d}$ is an ordering of the primitive idempotents of $A$.
\item $\{E^{*}_{i}\}_{i=0}^{d}$ is an ordering of the primitive idempotents of $A^{*}$.
\item ${\displaystyle{E^{*}_{i}AE^{*}_{j} =
\begin{cases}
0, & \text{if $\;|i-j|>1$;} \\
\neq 0, & \text{if $\;|i-j|=1$}
\end{cases}
}}
\qquad \qquad (0 \leq i,j\leq d).$
\item ${\displaystyle{E_{i}A^{*}E_{j} =
\begin{cases}
0, & \text{if $\;|i-j|>1$;} \\
\neq 0, & \text{if $\;|i-j|=1$}
\end{cases}
}}
\qquad \qquad (0 \leq i,j\leq d).$
\end{enumerate}
\end{definition}

\noindent Leonard systems and Leonard pairs are related as follows. Let $(A; \{E_{i}\}_{i=0}^{d}; A^{*}; \{E^{*}_{i}\}_{i=0}^{d})$ denote a Leonard system on $V$. For $0\leq i\leq d$, let $v_{i}$ denote a nonzero vector in $E_{i}V$. Then the sequence $\{v_{i}\}_{i=0}^{d}$ is a basis for $V$ which satisfies Definition \ref{def:lp}(ii). For $0\leq i\leq d$, let $v^{*}_{i}$ denote a nonzero vector in $E^{*}_{i}V$. Then the sequence $\{v^{*}_{i}\}_{i=0}^{d}$ is a basis for $V$ which satisfies Definition \ref{def:lp}(i). By these comments, the pair $A,A^{*}$ is a Leonard pair on $V$. Conversely, let $A,A^{*}$ denote a Leonard pair on $V$. By \cite[Lemma 1.3]{T:Leonard}, each of $A,A^{*}$ is multiplicity-free. Let $\{v_{i}\}_{i=0}^{d}$ denote a basis for $V$ which satisfies Definition \ref{def:lp}(ii). For $0\leq i\leq d$, the vector $v_{i}$ is an eigenvector for $A$; let $E_{i}$ denote the corresponding primitive idempotent. Let $\{v^{*}_{i}\}_{i=0}^{d}$ denote a basis for $V$ which satisfies Definition \ref{def:lp}(i). For $0\leq i\leq d$, the vector $v^{*}_{i}$ is an eigenvector for $A^{*}$; let $E^{*}_{i}$ denote the corresponding primitive idempotent. Then $(A; \{E_{i}\}_{i=0}^{d}; A^{*}; \{E^{*}_{i}\}_{i=0}^{d})$ is a Leonard system on $V$.

\medskip
\noindent We make some observations. Let $(A; \{E_{i}\}_{i=0}^{d}; A^{*}; \{E^{*}_{i}\}_{i=0}^{d})$ denote a Leonard system on $V$. For $0\leq i\leq d$, let $\theta_{i}$ (resp. $\theta^{*}_{i}$) denote the eigenvalue of $A$ (resp. $A^{*}$) associated with $E_{i}V$ (resp. $E^{*}_{i}V$). By construction, $\{\theta_{i}\}_{i=0}^{d}$ (resp. $\{\theta^{*}_{i}\}_{i=0}^{d}$) are mutually distinct and contained in $\fld$. It was shown in \cite[Lemma 12.7]{T:Leonard} that there exists $\beta\in\fld$ such that:
\begin{enumerate}
\item $\theta_{i-1}-\beta \theta_{i}+\theta_{i+1}$ is independent of $i$ for $1\leq i \leq d-1$;
\item $\theta^{*}_{i-1}-\beta \theta^{*}_{i}+\theta^{*}_{i+1}$ is independent of $i$ for $1\leq i \leq d-1$.
\end{enumerate}

\section{The antiautomorphism $\dagger$}

In this section, we discuss an antiautomorphism related to Leonard systems.

\begin{lemma} \label{lem:TDpower}
Let $A$ denote an irreducible tridiagonal matrix in $\mbox{\rm{Mat}}_{d+1}(\fld)$. Then the following {\rm (i)--(iii)} hold for $0\leq i,j\leq d$.
\begin{enumerate}
\item[\rm (i)] The entry $(A^{r})_{ij}=0$ if $r<|i-j|$\qquad \qquad $(0 \leq r \leq d)$.
\item[\rm (ii)] Suppose $i\leq j$. Then the entry $(A^{j-i})_{ij}=\prod_{h=i}^{j-1}A_{h,h+1}$. Moreover, $(A^{j-i})_{ij} \ne 0$.
\item[\rm (iii)] Suppose $i\geq j$. Then the entry $(A^{i-j})_{ij}=\prod_{h=j}^{i-1}A_{h+1,h}$. Moreover, $(A^{i-j})_{ij} \ne 0$.
\end{enumerate}
\end{lemma}
\noindent {\it Proof:} This follows from the definition of matrix multiplication and the meaning of irreducible tridiagonal. \hfill $\Box$ \\

\begin{assumption} \label{assum:E*AE*}
\rm
Let $\{E^{*}_{i}\}^{d}_{i=0}$ denote a system of mutually orthogonal idempotents in $\mathcal{A}$. Let $A$ denote an element of $\mathcal{A}$ such that
\begin{equation} \label{eq:E*AE*}
{\displaystyle{E^{*}_{i}AE^{*}_{j} =
\begin{cases}
0, & \text{if $\;|i-j|>1$;} \\
\neq 0, & \text{if $\;|i-j|=1$}
\end{cases}
}}
\qquad \qquad (0 \leq i,j\leq d).
\end{equation}
\end{assumption}

\begin{proposition} \label{lem:AE*Abasis}
With reference to Assumption \ref{assum:E*AE*}, the elements
\begin{equation}
A^{r}E^{*}_{0}A^{s} \qquad \qquad (0 \leq r,s\leq d) \label{eq:AE*Abasis}
\end{equation}
form a basis for the $\fld$-vector space $\mathcal{A}$.
\end{proposition}
\noindent {\it Proof:} We first show that the elements in the set (\ref{eq:AE*Abasis}) are linearly independent. To do this, we represent the elements in (\ref{eq:AE*Abasis}) by matrices. For $0\leq i\leq d$, let $v^{*}_{i}$ denote a nonzero vector in $E^{*}_{i}V$ and observe that $\{v^{*}_{i}\}_{i=0}^{d}$ is a basis for $V$. For $X \in \mathcal{A}$, let $X^{\flat}$ denote the matrix in $\mbox{Mat}_{d+1}(\fld)$ which represents $X$ with respect to the basis $\{v^{*}_{i}\}_{i=0}^{d}$. We observe that $\flat : {\mathcal A}\rightarrow \mbox{Mat}_{d+1}(\fld)$ is an isomorphism of $\fld$-algebras.  We abbreviate $B=A^{\flat}$, $F^{*}_{0}=E^{*\flat}_{0}$ and observe by (\ref{eq:E*AE*}) that $B$ is irreducible tridiagonal. For $0\leq r,s\leq d$, we show that the entries of $B^{r}F^{*}_{0}B^{s}$ satisfy
\begin{equation} \label{eq:AE*Aij1}
{\displaystyle{(B^{r}F^{*}_{0}B^{s})_{ij} =
\begin{cases}
0, & \text{if $i>r$ or $j>s$;} \\
\neq 0, & \text{if $i=r$ and $j=s$}
\end{cases}
}}
\qquad \qquad (0 \leq i,j\leq d).
\end{equation}
Because $\{E^{*}_{i}\}^{d}_{i=0}$ form a system of mutually orthogonal idempotents, $E^{*}_{0}v^{*}_{0}=v^{*}_{0}$ and $E^{*}_{0}v^{*}_{i}=0$ for $i\neq 0$. Therefore, the matrix $F^{*}_{0}$ has $(0,0)$-entry $1$ and all other entries $0$. So
\begin{equation}
(B^{r}F^{*}_{0}B^{s})_{ij}=(B^{r})_{i0}(B^{s})_{0j} \qquad \qquad (0 \leq i,j\leq d). \label{eq:AE*Aij2}
\end{equation}
Because $B$ is irreducible tridiagonal, Lemma \ref{lem:TDpower} applies. So, for $0\leq i\leq d$, the entry $(B^{r})_{i0}$ is zero if $i>r$ and nonzero if $i=r$. Similarly, for $0\leq j\leq d$, the entry $(B^{s})_{0j}$ is zero if $j>s$ and nonzero if $j=s$. Combining these facts with (\ref{eq:AE*Aij2}), we obtain (\ref{eq:AE*Aij1}), from which it follows that the elements in (\ref{eq:AE*Abasis}) are linearly independent. The number of elements in (\ref{eq:AE*Abasis}) is equal to $(d+1)^{2}$, which is the dimension of $\mathcal{A}$. Therefore, the elements in (\ref{eq:AE*Abasis}) form a basis for $\mathcal{A}$, as desired. \hfill $\Box$ \\

\begin{corollary} \label{cor:AE*0gen}
With reference to Assumption \ref{assum:E*AE*}, the elements $A$ and $E^{*}_{0}$ together generate $\mathcal{A}$.
\end{corollary}
\noindent {\it Proof:} This is an immediate consequence of Proposition \ref{lem:AE*Abasis}. \hfill $\Box$ \\

\noindent We recall the notion of an {\it antiautomorphism} of $\mathcal{A}$. Let $\gamma : \mathcal{A} \rightarrow \mathcal{A}$ denote any map. We call $\gamma$ an {\it antiautomorphism} of $\mathcal{A}$ whenever $\gamma$ is an isomorphism of $\fld$-vector spaces and $(XY)^{\gamma}=Y^{\gamma}X^{\gamma}$ for all $X, Y\in \mathcal{A}$.

\begin{lemma} \label{lem:dagexist}
With reference to Assumption \ref{assum:E*AE*}, there exists a unique antiautomorphism $\dagger$ of $\mathcal{A}$ such that $A^{\dagger}=A$ and $E^{*\dagger}_{0}=E^{*}_{0}$. Moreover, $E^{*\dagger}_{i}=E^{*}_{i}$ for $1 \leq i \leq d$ and $X^{\dagger \dagger}=X$ for all $X \in \mathcal{A}$. 
\end{lemma}
\noindent {\it Proof:} Concerning the existence of $\dagger$, we adopt the notation used in the proof of Proposition \ref{lem:AE*Abasis}. For $0 \leq i \leq d$, let $F^{*}_{i}=E^{*\flat}_{i}$ and note that $F^{*}_{i}$ is diagonal with $(i,i)$-entry $1$ and all other entries $0$. Recall that $B$ is irreducible tridiagonal. Let $D$ denote the diagonal matrix in $\mbox{Mat}_{d+1}(\fld)$ which has $(i,i)$-entry
\begin{equation}
D_{ii}= \frac{B_{01}B_{12}\cdots B_{i-1,i}}{B_{10}B_{21}\cdots B_{i,i-1}} \qquad \qquad (0 \leq i \leq d). \notag
\end{equation}
It is routine to verify $D^{-1}B^{t}D=B$, where $t$ denotes transpose. Fix an integer $i$ ($0\leq i\leq d$). Recall that $F^{*}_{i}$ is diagonal, so $F^{*t}_{i}=F^{*}_{i}$. Also, $D$ is diagonal, so $DF^{*}_{i}=F^{*}_{i}D$. From these comments, $D^{-1}F^{*t}_{i}D=F^{*}_{i}$. Define a map $\sigma :{\mbox{Mat}}_{d+1}(\fld) \rightarrow {\mbox{Mat}}_{d+1}(\fld)$ which satisfies $X^{\sigma}=D^{-1}X^{t}D$ for all $X \in {\mbox{Mat}}_{d+1}(\fld)$. We observe that $\sigma$ is an antiautomorphism of ${\mbox{Mat}}_{d+1}(\fld)$ such that $B^{\sigma}=B$ and $F^{*\sigma}_{i}=F^{*}_{i}$ for $0 \leq i \leq d$. We define the map $\dagger :\mathcal{A} \rightarrow \mathcal{A}$ to be the composition $\flat \sigma \flat^{-1}$. We observe that $\dagger$ is an antiautomorphism of $\mathcal{A}$ such that $A^{\dagger}=A$ and $E^{*\dagger}_{i}=E^{*}_{i}$ for $0 \leq i \leq d$. We have now shown that there exists an antiautomorphism $\dagger$ of $\mathcal{A}$ such that $A^{\dagger}=A$ and $E^{*\dagger}_{i}=E^{*}_{i}$ for $0 \leq i \leq d$. Our assertion about uniqueness follows from the fact that $A$ and $E^{*}_{0}$ together generate $\mathcal{A}$. The map $X\mapsto X^{\dagger \dagger}$ is an isomorphism of $\fld$-algebras from $\mathcal{A}$ to itself. This map is the identity since $A^{\dagger \dagger}=A$, $E^{*\dagger \dagger}_{0}=E^{*}_{0}$, and $\mathcal{A}$ is generated by $A$ and $E^{*}_{0}$.
\hfill $\Box$ \\

\noindent Up until now, we have been discussing the situation of Assumption \ref{assum:E*AE*}. We now modify this situation as follows.

\begin{assumption} \label{assum:A*}
\rm
Let $A$ and $\{E^{*}_{i}\}_{i=0}^{d}$ be as in Assumption \ref{assum:E*AE*}. Furthermore, assume that $A$ is multiplicity-free, with primitive idempotents $\{E_{i}\}_{i=0}^{d}$ and eigenvalues $\{\theta_{i}\}_{i=0}^{d}$. Additionally, let $\{\theta^{*}_{i}\}^{d}_{i=0}$ denote scalars in $\fld$ and let $A^{*}=\sum_{i=0}^{d}\theta^{*}_{i}E^{*}_{i}$. To avoid trivialities, assume that $d\geq 1$.
\end{assumption}

\begin{lemma} \label{lem:eistab}
With reference to Assumption \ref{assum:A*}, the antiautomorphism $\dagger$ from Lemma \ref{lem:dagexist} satisfies $A^{*\dagger}=A^{*}$ and $E_{i}^{\dagger}=E_{i}$ for $0 \leq i \leq d$.
\end{lemma}
\noindent {\it Proof:} By (\ref{eq:EpolyA}), $E_{i}$ is a polynomial in $A$ for $0 \leq i \leq d$. The result follows in view of Lemma \ref{lem:dagexist}. \hfill $\Box$ \\

\begin{lemma} \label{lem:undirected}
With reference to Assumption \ref{assum:A*} and for $0\leq i,j\leq d$, $E_{i}A^{*}E_{j}=0$ if and only if $E_{j}A^{*}E_{i}=0$.
\end{lemma}
\noindent {\it Proof:} Let $\dagger$ be the antiautomorphism from Lemma \ref{lem:dagexist}. Then $E_{i}A^{*}E_{j}=0$ if and only if $(E_{i}A^{*}E_{j})^{\dagger}=0$. Also, using Lemma \ref{lem:eistab}, $(E_{i}A^{*}E_{j})^{\dagger}=E_{j}^{\dagger}A^{* \dagger}E_{i}^{\dagger}=E_{j}A^{*}E_{i}$. The result follows. \hfill $\Box$ \\

\section{The graph $\Delta$} \label{sec:delta}

In the following discussion, a graph is understood to be finite and undirected, without loops or multiple edges.

\begin{definition} \label{def:delta}
\rm
With reference to Assumption \ref{assum:A*}, let $\Delta$ be the graph with vertex set $\{0, 1,\ldots ,d\}$ such that two vertices $i$ and $j$ are adjacent if and only if $i\neq j$ and $E_{i}A^{*}E_{j}\neq 0$. The graph $\Delta$ is well-defined in view of Lemma \ref{lem:undirected}.
\end{definition}

\begin{lemma} \label{lem:lspath}
With reference to Assumption \ref{assum:A*}, the following are equivalent:
\begin{enumerate}
\item[\rm (i)] the sequence $(A; \{E_{i}\}^{d}_{i=0}; A^{*}; \{E^{*}_{i}\}^{d}_{i=0})$ is a Leonard system;
\item[\rm (ii)] the graph $\Delta$ is a path such that vertices $i-1, i$ are adjacent for $1\leq i\leq d$.
\end{enumerate}
\end{lemma}
\noindent {\it Proof:} (i) $\Rightarrow$ (ii). This follows from condition (v) of Definition \ref{def:ls}.

\medskip
\noindent (ii) $\Rightarrow$ (i). We show that conditions (i)--(v) of Definition \ref{def:ls} are satisfied. Note that properties (ii) and (iv) of Definition \ref{def:ls} are satisfied by Assumption \ref{assum:A*}, while property (v) of Definition \ref{def:ls} is satisfied by construction. Concerning condition (i) of Definition \ref{def:ls}, we assume that $A$ is multiplicity-free. We now show that $A^{*}$ is multiplicity-free. Define a polynomial $m(\lambda)=\prod_{i=0}^{d}(\lambda -\theta^{*}_{i})$ and note that $m(A^{*})=0$ by Assumption \ref{assum:A*}. For $0\leq i\leq d$, let $v_{i}$ denote a nonzero vector in $E_{i}V$. Observe that $\{v_{i}\}_{i=0}^{d}$ is a basis for $V$. By construction, the matrix representing $A^{*}$ with respect to this basis is irreducible tridiagonal. The elements $\{A^{*i}\}_{i=0}^{d}$ are linearly independent by Lemma \ref{lem:TDpower}, so the minimal polynomial of $A^{*}$ has degree $d+1$. Therefore, the minimal polynomial of $A^{*}$ is precisely $m(\lambda)$. Because $A^{*}$ is diagonalizable, $m(\lambda)$ has distinct roots. It follows that $\{\theta^{*}_{i}\}_{i=0}^{d}$ are mutually distinct. Therefore, $A^{*}$ is multiplicity-free as desired. We have established condition (i) of Definition \ref{def:ls}. By Assumption \ref{assum:A*} and since $A^{*}$ is multiplicity-free, we see that $\{E^{*}_{i}\}_{i=0}^{d}$ is an ordering of the primitive idempotents of $A^{*}$. This gives property (iii) of Definition \ref{def:ls}. By these comments, $(A; \{E_{i}\}^{d}_{i=0}; A^{*}; \{E^{*}_{i}\}^{d}_{i=0})$ is a Leonard system. \hfill $\Box$ \\ 

\begin{definition} \label{def:Qpoly1}
\rm
With reference to Assumption \ref{assum:A*}, the given ordering $\{E_{i}\}_{i=0}^{d}$ of the primitive idempotents of $A$ is said to be {\it $Q$-polynomial}
whenever the equivalent conditions (i), (ii) hold in Lemma \ref{lem:lspath}.
\end{definition}

\begin{definition} \label{def:Qpoly2}
\rm
With reference to Assumption \ref{assum:A*}, let $(E,F)$ denote an ordered pair of distinct primitive idempotents for $A$. This pair will be called {\it $Q$-polynomial} whenever there exists a $Q$-polynomial ordering $\{E_{i}\}_{i=0}^{d}$ of the primitive idempotents of $A$ such that $E=E_{0}$ and $F=E_{1}$.
\end{definition}

\noindent The following is motivated by \cite[Definition 5.1]{Lang}.

\begin{definition} \label{def:tail}
\rm
With reference to Assumption \ref{assum:A*}, let $(E,F)=(E_{i},E_{j})$ denote an ordered pair of distinct primitive idempotents for $A$. This pair will be called a {\it tail} whenever the following occurs in $\Delta$:
\begin{enumerate}
\item[\rm (i)] $i$ is adjacent to no vertex in $\Delta$ besides $j$;
\item[\rm (ii)] $j$ is adjacent to at most one vertex in $\Delta$ besides $i$.
\end{enumerate}
\end{definition}

\begin{lemma} \label{lem:Qpolytail}
With reference to Assumption \ref{assum:A*}, let $(E,F)$ denote an ordered pair of distinct primitive idempotents for $A$. If $(E,F)$ is $Q$-polynomial, then $(E,F)$ is a tail.
\end{lemma}
\noindent {\it Proof:} Compare Definitions \ref{def:Qpoly1} and \ref{def:tail}. \hfill $\Box$ \\

\noindent For the rest of this section, we discuss the relationship between the connectivity of $\Delta$ and the subspaces of $V$ that are invariant under both $A$ and $A^{*}$.

\begin{lemma} \label{lem:invconnect1}
With reference to Assumption \ref{assum:A*}, fix a subspace $U\subseteq V$. Then $AU \subseteq U$ if and only if there exists a subset $S\subseteq \{0, 1,\ldots ,d\}$ such that $U=\sum_{h\in S}E_{h}V$. In this case, $S$ is uniquely determined by $U$.
\end{lemma}
\noindent {\it Proof:} First, assume there exists $S\subseteq \{0, 1,\ldots ,d\}$ such that $U=\sum_{h\in S}E_{h}V$. Then $AU\subseteq U$ since $AE_{i}=\theta_{i}E_{i}$ for $0\leq i\leq d$. Conversely, assume that $AU\subseteq U$. For $0\leq h\leq d$, we have $E_{h}U\subseteq U$ since $E_{h}$ is a polynomial in $A$. Therefore, $\sum_{h=0}^{d}E_{h}U\subseteq U$. Also, $U\subseteq \sum_{h=0}^{d}E_{h}U$ since $I=\sum_{h=0}^{d}E_{h}$. Therefore, $U=\sum_{h=0}^{d}E_{h}U$. Choose an integer $h$ ($0\leq h\leq d$). We have $E_{h}U\subseteq E_{h}V$ since $U\subseteq V$. The space $E_{h}V$ has dimension one, so $E_{h}U$ is either $0$ or $E_{h}V$. By these comments, there exists a subset $S\subseteq \{0, 1,\ldots ,d\}$ such that $U=\sum_{h\in S}E_{h}V$. It is clear that $S$ is uniquely determined by $U$. \hfill $\Box$ \\

\noindent We will use the following notation. For a subset $S\subseteq \{0, 1,\ldots ,d\}$, let $\overline{S}$ denote the complement of $S$ in $\{0, 1,\ldots ,d\}$.

\begin{proposition} \label{prop:invconnect2}
With reference to Assumption \ref{assum:A*}, fix a subset $S\subseteq \{0, 1,\ldots ,d\}$ and let $U=\sum_{h\in S}E_{h}V$. Then the following are equivalent: 
\begin{enumerate}
\item[\rm (i)] $A^{*}U \subseteq U$;
\item[\rm (ii)] the vertices $i,j$ are not adjacent in the graph $\Delta$ for all $i\in S$ and $j\in \overline{S}$.
\end{enumerate}
\end{proposition}
\noindent {\it Proof:} (i) $\Rightarrow$ (ii). Let $i\in S$ and $j\in \overline{S}$. Note that $E_{i}V\subseteq U$, so $E_{j}A^{*}E_{i}V\subseteq E_{j}A^{*}U\subseteq E_{j}U$ since $A^{*}U \subseteq U$. By assumption, $E_{j}U=E_{j}(\sum_{h\in S}E_{h}V)=0$ because $j\notin S$ and $E_{j}E_{h}=0$ for $j\neq h$. Thus, $E_{j}A^{*}E_{i}=0$, so $i$ and $j$ are not adjacent in $\Delta$.

\medskip
\noindent (ii) $\Rightarrow$ (i). It suffices to show that $A^{*}E_{i}V\subseteq U$ for $i\in S$. Let $i\in S$ be given. Using $\sum_{h=0}^{d}E_{h}=I$ and Definition \ref{def:delta}, we find $A^{*}E_{i}V=\sum_{h=0}^{d}E_{h}A^{*}E_{i}V=\sum_{h\in S}E_{h}A^{*}E_{i}V\subseteq \sum_{h\in S}E_{h}V=U$. The result follows. \hfill $\Box$ \\

\section{The main theorem}

The following is our main result.

\begin{theorem} \label{thm:main}
With reference to Assumption \ref{assum:A*}, let $(E,F)$ denote an ordered pair of distinct primitive idempotents for A. Then this pair is $Q$-polynomial if and only if the following {\rm (i)--(iii)} hold.
\begin{enumerate}
\item[\rm (i)] $(E,F)$ is a tail.
\item[\rm (ii)] There exists $\beta \in \fld$ such that $\theta^{*}_{i-1}-\beta \theta^{*}_{i}+\theta^{*}_{i+1}$ is independent of $i$ for $1\leq i \leq d-1$.
\item[\rm (iii)] $\theta^{*}_{0}\neq \theta^{*}_{i}$ for $1\leq i \leq d$.
\end{enumerate}
\end{theorem}
\noindent {\it Proof:} First, assume that $(E,F)$ is $Q$-polynomial. Condition (i) follows from Lemma \ref{lem:Qpolytail}. Conditions (ii) and (iii) follow from the last paragraph of Section \ref{sec:LS}.

\medskip
\noindent Conversely, assume that $(E,F)$ satisfies conditions (i)--(iii). We show that $(E,F)$ is $Q$-polynomial. To do this, we consider the graph $\Delta$ from Definition \ref{def:delta}. We begin by showing that $\Delta$ is connected. Suppose $\Delta$ is not connected. Then there exists a non-empty proper subset $S$ of $\{0,1,\ldots ,d\}$ such that $i$ and $j$ are not adjacent in $\Delta$ for all $i\in S$ and $j\in \overline{S}$. Let $U=\sum_{h\in S}E_{h}V$ and note that $U\neq 0$ and $U\neq V$. Observe that $AU\subseteq U$ by Lemma \ref{lem:invconnect1} and $A^{*}U\subseteq U$ by Proposition \ref{prop:invconnect2}. Using the equation $A^{*}=\sum_{i=0}^{d}\theta^{*}_{i}E^{*}_{i}$ and the fact that $\{E^{*}_{i}\}_{i=0}^{d}$ are mutually orthogonal idempotents,
\begin{equation}
E^{*}_{0}=\prod_{j=1}^{d}\frac{A^{*}-\theta^{*}_{j}I}{\theta^{*}_{0}-\theta^{*}_{j}}. \label{eq:E*0polyA*}
\end{equation}
Note that the denominator is nonzero by condition (iii). By (\ref{eq:E*0polyA*}) and since $A^{*}U\subseteq U$, we find that $E^{*}_{0}U\subseteq U$. By Corollary \ref{cor:AE*0gen}, $A$ and $E^{*}_{0}$ generate $\mathcal{A}$. Therefore, $\mathcal{A}U\subseteq U$. Recall that $V$ is irreducible as an $\mathcal{A}$-module, so either $U=0$ or $U=V$. This is a contradiction, so $\Delta$ is connected.

\medskip
\noindent Relabeling the primitive idempotents of $A$ as necessary, we may assume without loss of generality that $E_{0}=E$ and $E_{1}=F$. Because $(E,F)$  is a tail and $\Delta$ is connected, vertex $0$ is adjacent to vertex $1$ and no other vertices. Similarly, vertex $1$ is adjacent to vertex $0$ and at most one other vertex. We now show that $\Delta$ is a path.

\medskip
\noindent First, let $\gamma^{*}$ be the common value of $\theta^{*}_{i-1}-\beta \theta^{*}_{i}+\theta^{*}_{i+1}$ for $1\leq i \leq d-1$. We claim that the expression
\begin{equation}
\theta ^{*2}_{i-1}-\beta \theta ^{*}_{i-1}\theta ^{*}_{i}+\theta ^{*2}_{i}-\gamma ^{*}(\theta ^{*}_{i-1}+\theta ^{*}_{i}) \label{eq:delta*}
\end{equation}
is independent of $i$ for $1\leq i\leq d$. Let $p_{i}$ denote expression (\ref{eq:delta*}). Observe that, for $1\leq i\leq d-1$,
\begin{equation}
p_{i}-p_{i+1}=(\theta^{*}_{i-1}-\theta^{*}_{i+1})(\theta^{*}_{i-1}-\beta\theta^{*}_{i}+\theta^{*}_{i+1}-\gamma^{*}), \notag
\end{equation}
which therefore equals $0$. Consequently, $p_{i}$ is independent of $i$ for $1\leq i\leq d$. The claim is now proved. Let $\delta^{*}$ denote the common value of (\ref{eq:delta*}) for $1\leq i\leq d$. We now show that
\begin{equation}
0=[A^{*},A^{*2}A-\beta A^{*}AA^{*}+AA^{*2}-\gamma ^{*}(AA^{*}+A^{*}A)-\delta ^{*}A], \label{eq:TD2}
\end{equation}
where $[x,y]=xy-yx$.

\medskip
\noindent Let $C$ denote the expression on the right-hand side of (\ref{eq:TD2}). Using $I=\sum_{i=0}^{d}E^{*}_{i}$, we obtain
\begin{align}
C &= (E^{*}_{0}+E^{*}_{1}+\cdots +E^{*}_{d})C(E^{*}_{0}+E^{*}_{1}+\cdots +E^{*}_{d}) \notag \\
  &= \sum_{i=0}^{d}\sum_{j=0}^{d} E^{*}_{i}CE^{*}_{j}. \notag
\end{align}
To show that $C=0$, it suffices to show that $E^{*}_{i}CE^{*}_{j}=0$ for $0\leq i,j\leq d$. Let $i$ and $j$ be given. Recall that $E^{*}_{i}A^{*}=\theta^{*}_{i}E^{*}_{i}$ and $A^{*}E^{*}_{j}=\theta^{*}_{j}E^{*}_{j}$.
Thus,
\begin{equation}
E^{*}_{i}CE^{*}_{j}=(E^{*}_{i}AE^{*}_{j})P(\theta^{*}_{i},\theta^{*}_{j})(\theta^{*}_{i}-\theta^{*}_{j}), \notag
\end{equation}
where
\begin{equation}
P(\lambda,\mu)=\lambda^{2}-\beta \lambda \mu+\mu^{2}-\gamma^{*}(\lambda +\mu)-\delta^{*}. \notag
\end{equation}
 
\medskip
\noindent If $|i-j|>1$, then $E^{*}_{i}AE^{*}_{j}=0$ by Assumption \ref{assum:A*}. If $|i-j|=1$, then $P(\theta^{*}_{i},\theta^{*}_{j})=0$. If $i=j$ then $\theta^{*}_{i}-\theta^{*}_{j}=0$. Therefore, $E^{*}_{i}CE^{*}_{j}=0$ in all cases, so $C=0$. We have now shown (\ref{eq:TD2}).

\medskip
\noindent Suppose we are given vertices $i$ and $j$ in $\Delta$ at $\partial (i,j)=3$, where $\partial$ denotes path-length distance. Further, suppose there exists a unique path of length $3$ connecting $i$ and $j$. Denoting this path by $(i,r,s,j)$, we show 
\begin{equation}
\theta_{i}-(\beta +1)\theta_{r}+(\beta +1)\theta_{s}-\theta_{j}=0. \label{eq:irsj}
\end{equation}
To show (\ref{eq:irsj}), expand the right-hand side of (\ref{eq:TD2}) to get
\begin{align}
0 = &A^{*3}A-(\beta+1)A^{*2}AA^{*}+(\beta+1)A^{*}AA^{*2}-AA^{*3} \notag \\
    &-\gamma^{*}(A^{*2}A-AA^{*2})-\delta^{*}(A^{*}A-AA^{*}). \notag
\end{align}
In the above equation, multiply each term on the left by $E_{i}$ and on the right by $E_{j}$, and simplify. To illustrate, we now simplify the first term. Using $AE_{j}=\theta_{j}E_{j}$, we find that $E_{i}A^{*3}AE_{j}=\theta_{j}E_{i}A^{*3}E_{j}$. Using Lemma \ref{lem:EsumI},
\begin{align}
E_{i}A^{*3}E_{j} &= E_{i}A^{*}\left( \sum_{h=0}^{d}E_{h} \right) A^{*} \left( \sum_{l=0}^{d}E_{l} \right) A^{*}E_{j} \notag \\
                 &= E_{i}A^{*}E_{r}A^{*}E_{s}A^{*}E_{j}. \notag
\end{align}
Therefore, 
\begin{equation}
E_{i}A^{*3}AE_{j}=\theta_{j}E_{i}A^{*}E_{r}A^{*}E_{s}A^{*}E_{j}. \notag
\end{equation}
Simplifying the other terms in a similar fashion yields
\begin{align}
E_{i}A^{*2}AA^{*}E_{j} &= \theta_{s}E_{i}A^{*}E_{r}A^{*}E_{s}A^{*}E_{j}, \notag \\
E_{i}A^{*}AA^{*2}E_{j} &= \theta_{r}E_{i}A^{*}E_{r}A^{*}E_{s}A^{*}E_{j}, \notag \\
     E_{i}AA^{*3}E_{j} &= \theta_{i}E_{i}A^{*}E_{r}A^{*}E_{s}A^{*}E_{j}, \notag
\end{align}
\begin{align}
E_{i}A^{*2}AE_{j} &= 0, & E_{i}AA^{*2}E_{j}=0, \notag \\
 E_{i}A^{*}AE_{j} &= 0, &  E_{i}AA^{*}E_{j}=0. \notag
\end{align}
By the above comments, we get 
\begin{equation}
0 = \big( \theta_{i}-(\beta+1)\theta_{r}+(\beta+1)\theta_{s}-\theta_{j}\big) E_{i}A^{*}E_{r}A^{*}E_{s}A^{*}E_{j}. \label{eq:TD2result}
\end{equation}
Since $s$ and $j$ are adjacent, $E_{s}A^{*}E_{j}\neq 0$. Therefore, $E_{s}A^{*}E_{j}V$ is a nonzero subspace of the one-dimensional space $E_{s}V$, so it follows that $E_{s}A^{*}E_{j}V=E_{s}V$. Similarly, $E_{r}A^{*}E_{s}V=E_{r}V$ and $E_{i}A^{*}E_{r}V=E_{i}V$, so $E_{i}A^{*}E_{r}A^{*}E_{s}A^{*}E_{j}V=E_{i}V$. Therefore, $E_{i}A^{*}E_{r}A^{*}E_{s}A^{*}E_{j}\neq 0$. This and (\ref{eq:TD2result}) imply (\ref{eq:irsj}).

\medskip
\noindent We can now easily show that $\Delta$ is a path. To this end, we show that every vertex in $\Delta$ is adjacent to at most two other vertices. Suppose there exists a vertex $i$ in $\Delta$ that is adjacent to at least three other vertices. Choose the $i$ such that $\partial (0,i)$ is minimum. Without loss of generality, assume that the vertices of $\Delta$ are labelled such that $\partial (0,i)=i$ and $(0,1,\ldots,i)$ is a path. By construction, $i\geq 2$. By assumption, there exist distinct vertices $j$ and $j'$, each at least $i+1$, that are both adjacent to $i$. Note that $\partial (i-2,j)=3$ and that $(i-2,i-1,i,j)$ is the unique path of length $3$ connecting $i-2$ and $j$. Therefore, by (\ref{eq:irsj}), 
\begin{equation}
\theta_{i-2}-(\beta +1)\theta_{i-1}+(\beta +1)\theta_{i}-\theta_{j}=0. \label{eq:thetaj}
\end{equation}
Replacing $j$ by $j'$ in the above argument, we obtain 
\begin{equation}
\theta_{i-2}-(\beta +1)\theta_{i-1}+(\beta +1)\theta_{i}-\theta_{j'}=0. \label{eq:thetaj'}
\end{equation}
Comparing (\ref{eq:thetaj}) to (\ref{eq:thetaj'}), we find $\theta _{j}=\theta _{j'}$. Recall that $\{\theta_{h}\}_{h=0}^{d}$ are mutually distinct, so $j=j'$. This is a contradiction and we have now shown that $\Delta$ is a path.

\medskip
\noindent The ordering of primitive idempotents $E_{0},E_{1},\ldots$ induced by the path is $Q$-polynomial by Definition \ref{def:Qpoly1}. Now the pair $(E,F)=(E_{0},E_{1})$ is $Q$-polynomial in view of Definition \ref{def:Qpoly2}.  \hfill $\Box$ \\

\section{Acknowledgment}

This paper was written while the author was a graduate student at the University of Wisconsin-Madison. The author would like to thank his advisor, Paul Terwilliger, for offering many valuable ideas and suggestions.

\bigskip

\noindent Edward Hanson \hfil\break
\noindent Department of Mathematics \hfil\break
\noindent University of Wisconsin \hfil\break
\noindent 480 Lincoln Drive \hfil\break
\noindent Madison, WI 53706-1388 USA \hfil\break
\noindent email: {\tt hanson@math.wisc.edu }\hfil\break

\end{document}